\documentclass{amsart}                                                          
\usepackage{amsxtra,amssymb,multicol,graphicx}                                  
\title{WHAT IS a Perverse Sheaf?}                                               
\author{Mark Andrea  de Cataldo and Luca Migliorini}                                     
\begin{document}                                                                
\rightline{\Huge \bf WHAT IS a Perverse Sheaf?}                                 
\vskip1cm                                                                       
\rightline{\it \huge Mark  Andrea de Cataldo and Luca Migliorini}                        
\vskip1in

\begin{multicols}{2}                                                          
Manifolds are  obtained  by glueing  open subsets
 of Euclidean space. Differential forms,  vector fields, etc.
 are  defined locally and  then  glued    to yield a 
 global object.                                                                         
The notion of sheaf 
embodies  the idea of glueing.
Sheaves come in many flavors:
 sheaves  of differential forms, of
vector fields, of differential operators,  
constant and locally constant sheaves, etc.
A locally constant sheaf ({\em local system}) on a space $X$  is determined 
by its monodromy,  i.e. by a representation of the fundamental group
$\pi_1(X,x)$
into the group of  automorphisms of  the fiber  at $x\in X$:
the sheaf of orientations on the M\"obius strip assigns $-{\rm Id}$ to
the   generators of the fundamental group ${\mathbb Z}$.
A sheaf, or even a map of sheaves, can be glued back together
from its local data:
 exterior derivation can be viewed  as a map between sheaves of differential forms;
the glueing is possible because exterior derivation is independent
   of the choice of local coordinates. 
 
 The theory of sheaves is made more complete  by considering {\em
 complexes of sheaves}. A complex  of sheaves  $K$ is  a collection
 of sheaves $\{K^i\}_{i \in {\mathbb Z}}$ and maps $d^i: K^i \to K^{i+1}$ subject 
 to $d^2 =0$. The $i$-th {\em cohomology sheaf}
 ${\mathcal H}^i(K)$ is  $ \mbox{ker} \,  d^i/\mbox{im} \,d^{i-1}$.
 The (sheafified) de Rham complex  ${\mathcal E}$
is the  complex with entries  the sheaves 
${\mathcal E}^i$ of differential $i$-forms and with
differentials $d: {\mathcal E}^i \to {\mathcal E}^{i+1}$ given by the  exterior derivation of differential forms.  By the Poincar\'e lemma,
the cohomology sheaves 
are all zero, except for  ${\mathcal H}^0 \simeq {\mathbb C},$ the constant sheaf.

   The de Rham theorem, stating that the cohomology of the constant sheaf equals
closed forms modulo exact ones, points to the fact that ${\mathbb C}$ and ${\mathcal E}$
are cohomologically indistinguishable
from each other, even at the local level. 
The need to identify    two complexes   
containing the same cohomological information via an  isomorphism leads to the notion
of  {\em derived category} (\cite{illu}): the objects are
complexes and  the arrows are designed to achieve
the desired identifications.
The   inclusion of complexes ${\mathbb C} \subseteq
{\mathcal E}$  is promoted by decree 
 to the rank of   isomorphism in the derived category
 because it induces an isomorphism at the level of cohomology sheaves.

 While  the derived category brings in a  
thick layer of abstraction, 
it   extends
the reach and  flexibility  of the theory.   
 One  defines the  cohomology
groups  of  a complex and  
extends to complexes of sheaves 
the ordinary
operations of  algebraic topology:  pull-backs,
push-forwards, cup  and cap  products, etc.
There is also 
a general form of duality  for complexes  (\cite{illu}) generalizing
 classical Poincar\'e duality.

Perverse sheaves live on spaces with singularities:
analytic spaces, algebraic varieties, PL spaces, pseudo-manifolds, etc.
For ease of exposition,
we  limit ourselves  to
sheaves of vector spaces on complex algebraic varieties
and to  perverse sheaves with respect to what is called {\em middle perversity}.
In order to avoid dealing with
pathologies such as sheaves supported on the Cantor set, 
one  imposes a technical condition called {\em constructibility}. 
Let us just say
that the category $D_X$  of
{\em bounded constructible complexes} of sheaves on $X$ sits 
in  the derived category
and is stable under the various topological  operations mentioned above.
If $K$ is in $D_X,$  only finitely many of its cohomology sheaves 
are non-zero and, for every $i$, the set $\mbox{ supp}\,{\mathcal H}^i(K)$,
the closure of the set
of points at which the stalk is non-zero, 
is an algebraic subvariety.

A {\em perverse sheaf} on  $X$  is a bounded  constructible complex $P \in D_X$ such that
the following holds for $K=P$ and for its dual $P^{\vee}$: 
\begin{equation}
\label{req}
\dim_{\mathbb C} {\rm supp}\,{\mathcal H}^{-i} (K) \leq i, \;\;\; \forall i \in {\mathbb Z}.
\end{equation} 
A map of perverse sheaves is an arrow in $D_X$.

 The term ``sheaf" stems from  the fact that,  just like in the case of ordinary   sheaves,
(maps of) perverse sheaves can be glued; as to ``perverse," see below.  
The theory of perverse sheaves
has its roots  in the two notions of intersection cohomology and
 of ${\mathcal D}$-module. As we see below, perverse  sheaves and
${\mathcal D}$-modules are related by the Riemann-Hilbert correspondence.

It is time for  examples.
 If $X$ is nonsingular, then ${\mathbb C}_X[\dim{X}]$,
  i.e. the constant sheaf in degree $-\dim_{\mathbb C}{X}$,
is self-dual and perverse. 
If $Y \subseteq X$ is a nonsingular closed subvariety, then
${\mathbb C}_Y[\dim{Y}]$, viewed as a complex on $X$, is a  perverse sheaf on $X$.
If $X$ is singular,
then ${\mathbb C}_X [\dim{X}]$ is usually not a perverse  sheaf.
On the other hand, the   intersection cohomology complex (see below) is a perverse sheaf,
regardless of the singularities of $X$.
The extension of two perverse sheaves is  a perverse sheaf.
The following example  can serve as a test case for the
first definitions in the theory of ${\mathcal D}$-modules. 
Let $X={\mathbb C}$ be the complex line with origin $o \in X$,
let
$z$ be the  standard holomorphic coordinate, let ${\mathcal O}_X$ be the sheaf
of holomorphic functions on $X$, let $a$ be a complex number and let $D$
 be the  differential operator $D:f\mapsto zf'- af$. The complex 
 $P_a$
 \begin{equation}
 \label{exps}
   0 \longrightarrow P^{-1}_a
  :={\mathcal O}_X \stackrel{D}\longrightarrow P^0_a:= {\mathcal O}_X 
 \longrightarrow 0\end{equation}
 is perverse. If $a \in {\mathbb Z}^{\geq 0}$, then 
 ${\mathcal H}^{-1} (P_a) = {\mathbb C}_X$ and ${\mathcal H}^{0}(P_0)= {\mathbb C}_{o}$.
If $a \in {\mathbb Z}^{< 0}$, then 
 ${\mathcal H}^{-1} (P_a)$ is the extension by zero at $o$ of the sheaf ${\mathbb C}_{X \setminus o}$
 and  ${\mathcal H}^{0}(P_a)= 0$.
 If $a \notin {\mathbb Z}$, then
  ${\mathcal H}^{-1} (P_a)$ is the extension 
by zero at $o$ of the local system on $X\setminus o$ associated with the branches of the  multi-valued function
$z^a$ and ${\mathcal H}^{0}(P_a)= 0$.   In each case, the associated  monodromy 
sends 
the  positive  generator of $\pi_1 (X \setminus o, 1)$ to
$e^{2\pi i a}$. The dual of $P_{a}$ is $P_{-a}$ (this fits well
with the notions of adjoint differential equation and of duality for ${\mathcal D}$-modules). Every $P_a$ is
the extension of the perverse sheaf ${\mathcal H}^{0}(P_a)[0]$  by the perverse sheaf ${\mathcal H}^{-1}(P_a)[1]$.
The extension is trivial (direct sum) if and only if $a\notin  {\mathbb Z}^{> 0}.$

A local system on a nonsingular variety can be turned 
into a perverse sheaf   
by viewing  it as 
a complex with   a single  entry in the appropriate  degree.
On the other hand, a perverse sheaf restricts to a local system
on some dense open subvariety.
We want to make sense of the following slogan:
{\em perverse sheaves are the
 singular version of  local systems.} In order to do so, we 
 discuss the two  widely different   ideas that led to the
 birth of perverse sheaves about thirty years ago: {\em the generalized
 Riemann-Hilbert correspondence} (RH) and {\em intersection cohomology} (IH)
 (\cite{kle}).

{\bf  (RH)}
Hilbert's $21^{st}$ problem is concerned with Fuchs-type
differential equations  on a  punctured  Riemann surface $\Sigma$. 
As one circuits the punctures, the solutions are transformed:
the sheaf of solutions is   a local system on $\Sigma$
(see (\ref{exps})).

The $21^{st}$ problem
asked whether any local system arises in this way (it  essentially does).
The sheafification
of  linear  partial differential equations on a manifold gives rise to the notion of
{\em ${\mathcal D}$-module}.
A  {\em regular holonomic} ${\mathcal D}$-module
on a complex manifold $M$ is the generalization of 
the Fuchs-type equations on $\Sigma$.
The sheaf of solutions is now  replaced by a {\em complex of solutions}
which, remarkably, belongs to $D_M$.  
In (\ref{exps}), the complex of solutions is
$P_a$,
the sheaf of solutions to $D(f)=0$  is ${\mathcal H}^{-1}(P_a)$ 
 and ${\mathcal H}^0(P_a)$ is related to
the (non)solvability of $D(f)=g$.
Let ${\mathcal D}^b_{r,h}(M)$ 
be the bounded derived category of ${\mathcal D}$-modules on $M$
 with regular holonomic cohomology.
 RH  states that the assignment of the  (dual to the) complex of solutions
 yields an equivalence of categories ${\mathcal D}^b_{r,h}(M) \simeq D_M$.
Perverse sheaves  enter the center of the stage:
they  correspond via RH  to
 regular holonomic ${\mathcal D}$-modules (viewed as complexes concentrated in degree 
 zero)
 
  In agreement with the slogan mentioned above, 
 the category of perverse sheaves
 shares the following formal properties with the category of local systems:
 it  is Abelian
(kernels, cokernels, images and coimages exist and
the  coimage is isomorphic to the image), stable under duality, Noetherian
 (the ascending chain condition holds) and
 Artinian (the descending chain condition
 holds), i.e. every perverse sheaf 
 is a finite iterated extension of {\em simple} (no subobjects) perverse sheaves. 
In our example, the perverse sheaves (\ref{exps}) are  simple if and only if
$a \in {\mathbb C} \setminus {\mathbb Z}.$

 What are the simple perverse sheaves? Intersection cohomology provides the answer.

{\bf  (IH)}
 The {\em intersection cohomology groups} of a singular variety $X$
 with coefficients in a local system are a topological invariant
 of the variety. 
 They coincide with ordinary cohomology when $X$ is nonsingular and the coefficients
 are constant.  
 These groups were originally  defined and studied
 using the  theory of geometric chains 
in order to study the failure, due to the presence of singularities, of  Poincar\'e duality for ordinary homology, and to  put a remedy to it by considering 
the homology theory arising by considering only
chains  that intersect the singular set in a controlled way.
  In this context, certain sequences of integers,   called   {\em perversities},
were  introduced to give a measure of
how a chain intersects the singular set, whence 
the origin of the term ``perverse."
The intersection cohomology groups thus defined 
satisfy the conclusions  of Poincar\'e duality and of  the Lefschetz hyperplane theorem.  

On the other hand, the intersection cohomology  groups  can also be exhibited
 as the cohomology groups of certain 
 complexes in $D_X$:
 the {\em intersection complexes}  of $X$  with   coefficients in the local system.
It is  a remarkable twist in the plot of this story, that the simple perverse sheaves
 are precisely the intersection complexes of the  irreducible
subvarieties   of $X$ with  coefficients given by simple local systems!

 We  are  now in a position to  clarify the  earlier slogan.
 A  local system  $L$ on a nonsingular subvariety $Z \subseteq M$ 
gives rise  to a regular holonomic ${\mathcal D}$-module
 supported over the closure $\overline{Z}$.
The same  $L$
gives rise to  the intersection complex of $\overline{Z}$ with coefficients in $L$.
Both objects extend $L$ from $Z$ to $\overline{Z}$
across the singularities $\overline{Z} \setminus Z$.
By RH: the intersection complex
is precisely the  complex of solutions of the  ${\mathcal D}$-module.

A pivotal role  in the applications of the theory of perverse sheaves
 is played by the {\em decomposition theorem}: 
 let $f: X \to Y$ be a proper map of varieties;
then the  intersection cohomology groups of $X$  with coefficients in a simple local system
   are isomorphic
to  the direct sum of a collection of   intersection cohomology
groups of irreducible subvarieties of $Y,$  with coefficients in simple local systems.
For example, if $f:X \to Y$ is a resolution of the singularities of 
$Y$, then  the intersection cohomology groups
of $Y$ are a direct summand of
the ordinary cohomology groups of  $X$.
This ``as-simple-as-possible" splitting behavior is 
 the deepest known fact concerning the 
homology of complex algebraic varieties and maps.
It   fails in complex analytic  and in real algebraic geometry. 
The decomposition
of the  intersection cohomology groups  of $X$ is a reflection in cohomology of a 
finer decomposition of complexes  in 
$D_Y$.
The original proof of the decomposition theorem
uses algebraic geometry over finite fields 
(perverse sheaves make perfect sense in this context). For a
discussion
of some of the  proofs see \cite{bams}.

One striking application of this circle of ideas is the fact that the intersection cohomology groups
of    projective varieties enjoy the same classical properties
of the  cohomology  groups of projective manifolds:  the Hodge $(p,q)$-decomposition theorem, 
the hard Lefschetz theorem, and the Hodge-Riemann bilinear relations. This, of course,
in addition  to Poincar\'e duality  and to  the Lefschetz hyperplane theorem  mentioned above.

The applications of the theory of  perverse sheaves range
from geometry to combinatorics to algebraic analysis.
The most dramatic ones are in the realm of
representation theory, where their introduction has  led to a
truly spectacular revolution: 
 proofs of the Kazhdan-Lusztig conjecture,  of the geometrization of the Satake isomorphism
and, recently,  of the fundamental lemma in the Langlands' program
(see the survey \cite{bams}).

\end{multicols} 


\begin{thebibliography}{99}
  
   \bibitem{bams}{ M.A. de Cataldo, L. Migliorini, ``The decomposition theorem,
   perverse sheaves and the topology of algebraic maps,"  a survey, 
   Bulletin of the A.M.S., Vol. 46, n.4, (2009), 535-633}

\bibitem{illu} {L. Illusie, ``Cat\'egories deriv\'ees et dualit\'e, travaux de J.L.Verdier'' Enseign.Math. {\bf (2) 36} (1990), 369-391.}

\bibitem{kle} {S. Kleiman, ``The development of intersection homology theory," Pure and Applied Mathematics Quarterly
vol. 3, n. 1, Special issue in honor of Robert MacPherson, 225-282, 2007.}


  
    \end{thebibliography}
\end{document}